\documentstyle[11pt]{article}

\input amssym.def
\input amssym

\def \C{{\Bbb C}}
\def \CC{{\cal C}}
\def \CF{{\cal F}}
\def \CP{{\cal P}}
\def \Fix{{\rm Fix}}
\def \Ga{\Gamma}
\def \ind{{\rm ind}}
\def \mod{{\rm mod \hspace{1pt}}}
\def \n{{\frak n}}
\def \ph{\varphi}
\def \prf{{\bf Proof: }}
\def \qed{\hfill $\Box$

$ $

}
\def \ra{\rightarrow}
\def \R{{\Bbb R}}
\def \Rep{{\rm Rep}}
\def \sign{{\rm sign}}
\def \tr{{\rm tr \hspace{1pt}}}
\def \Z{{\Bbb Z}}

\begin{document}

\title{A Lefschetz formula for flows}

\author{Anton Deitmar\\ {\small Math. Inst. d. Univ., INF 288,
69126 Heidelberg, Germany}}
\date{}
\maketitle

\pagestyle{myheadings}
\markright{A LEFSCHETZ FORMULA FOR FLOWS}

\tableofcontents

\newcommand{\rez}[1]{\frac{1}{#1}}
\newcommand{\der}[1]{\frac{\partial}{\partial #1}}
\newcommand{\binom}[2]{\left( \begin{array}{c}#1\\#2\end{array}\right)}

\newcounter{lemma}
\newcounter{corollary}
\newcounter{proposition}
\newcounter{theorem}

\newtheorem{conjecture}{\stepcounter{lemma} \stepcounter{corollary}
	\stepcounter{proposition}\stepcounter{theorem}
	Conjecture}[section]
\newtheorem{lemma}{\stepcounter{conjecture}\stepcounter{corollary}
	\stepcounter{proposition}\stepcounter{theorem}
	Lemma}[section]
\newtheorem{corollary}{\stepcounter{conjecture}\stepcounter{lemma}
	\stepcounter{proposition}\stepcounter{theorem}
	Corollary}[section]
\newtheorem{proposition}{\stepcounter{conjecture}\stepcounter{lemma}
	\stepcounter{corollary}\stepcounter{theorem}
	Proposition}[section]
\newtheorem{theorem}{\stepcounter{conjecture} \stepcounter{lemma}
	\stepcounter{corollary}\stepcounter{proposition}	Theorem}[section]

\begin{center}
{\bf Introduction}
\end{center}
Let M be a smooth compact manifold and $f:M\ra M$ a diffeomorphism with
nondegenerate fixed points.
Suppose we are given a flat vector bundle $E\ra M$ and that $f$ lifts linearly
to $E$.
Then $f$ acts by pullback on the $E$-valued differential forms and on the de
Rham cohomology $H^.(E)$ with coefficients in $E$.
The theorem of Atiyah-Bott-Lefschetz now tells us that
\begin{eqnarray}\label{ABL}
\sum_{a\in \Fix (f)} \ind_f(a) \tr (f\mid E_a) = \sum_{j=0}^{\dim M} (-1)^j
\tr(f^*\mid H^j(E)).
\end{eqnarray}
We want to formulate an analogous statement for a flow $\Phi$ instead of a
diffeomorphism.

By looking at the local side of (\ref{ABL}) the first idea would be to
consider the fixed points of the flow and to assume they are nondegenerate.
This way we get the Hopf formula expressing the Euler characteristic on the
global side.

More subtle information is supposed to be obtained by considering the closed
orbits of the flow.
The non-degeneracy of the latter translates to a hyperbolicity condition of the
flow.
Assume the flow lifts to the flat bundle $E$ then on each closed orbit $c$
with length $l(c)$ and on any point $m$ on $c$ we have a trace
$\tr(\Phi_{l(c)}\mid E_y)$ which is independent of $y$.
The index is at the first glance to be replaced by the Lefschetz index
$\ind_L(c)$ of the Poincar\'e map around the orbit, but it turns out that one
has to take into account the multiplicity $\mu(c)$ so we define the Fuller
index to be $\ind_F(c) := \ind_L(c)/\mu(c)$.
Our local side would then be the sum over all closed orbits of
$\ind_F(c)\tr\ph(c)$. Unfortunately this sum does not converge.
So we have to replace it by a zeta-regularized version
$\hat{\sum_{c}} \ind_F(c)\tr(\Phi_{l(c)}\mid E)$.

D. Fried \cite{Fr-tors} has shown that in case of geodesic flows on
hyperbolic spaces we have an identity
$$
-\hat{\sum_{c}} \ind_F(c)\tr(\Phi_{l(c)}\mid E) = \log \tau(E)\ \ \ \mod\ 2\pi
i\Z,
$$
where $\tau(E)$ is the analytic torsion of the bundle $E$.
This formula is an analogue of the Hopf formula but not of the
Atiyah-Bott-Lefschetz formula since on the global side the presence of the flow
is not visible.

In this paper we present an analogue of (\ref{ABL}).
Let $H$ be the generator of the flow.
We define a twisted version $H^.(\CC_\ph)$ of the tangential cohomology of the
stable foliation.
Then our Lefschetz formula is
$$
 \sum_{p=0}^{\dim Y -1} (-1)^{p} \log \det (H|H^P(\CC_\ph))=
 -\hat{\sum_{c}} \ind_F(c)\tr(\Phi_{l(c)}\mid E) \ \ \ \mod\ 2\pi i\Z,
$$
where det means the regularized determinant in the sense of D. Ray and I.
Singer \cite{RS-antors}.

\section{The main theorem}\label{sect_main_theo}
Let $X$ denote a compact hyperbolic manifold of dimension $d$, i.e. $X$ is a
Riemannian manifold of constant negative curvature which we may normalize to be
$-1$.
We will further restrict to the case that the dimension $d$ of $X$ is odd.

Now suppose $\Phi$ is the geodesic flow of $X$.
Then $\Phi$ acts on the sphere bundle $Y=SX$ of $X$.
It is known that $\Phi$ is Anosov, i.e. there is a smooth splitting of the
tangent bundle of $Y$ as
$$
TY = T_0 \oplus T_u \oplus T_s,
$$
where $T_0$ is spanned by the flow, the flow acts contractingly on the {\bf
stable part} $T_s$,  as the time tends to $+\infty$ and contractingly on the
{\bf unstable part} $T_u$,  as the time tends to $-\infty$.
The bundles $T_u$ and $T_s$ are integrable but their sum $T_s\oplus T_u$ is
not.
For a proof of the Anosov property see for example \cite{angeo}

Let the {\bf periodic set} $\CP$ of $\Phi$ be the subset of (space $\times$
time) $=Y\times (0,\infty)$ consisting of points $(y,t)$ such that $\Phi_t y =
y$.
For such a point $(y,t)$ let $c= \{ (\Phi_s y,s) | 0\leq s \leq t \}$ be the
underlying {\bf closed orbit} of $\Phi$, then $l(c) := t$ will be called the
{\bf length} or {\bf period} of $c$.

For any closed orbit $c$ of period $l(c)$ consider the linear Poincar\'e map
$$
P_{c,y} := D\Phi_{l(c)}|_{T_{s,y} \oplus T_{u,y}}
$$
for some point $(y,t)$ in $c$.
We are only interested in the index of this map, so the dependence on the point
will vanish.
Let
$$
\ind_L (c) := \sign \det (1-P_{c,y})
$$
be the {\bf Lefschetz index} of the orbit $c$.

For a closed orbit $c$ let $m(c)$ be its {\bf multiplicity}, i.e. the largest
natural number $m$ such that $t/m$ is a period of $y$, where $(y,t)\in c$.
Now the {\bf Fuller index} of $c$ is by definition
$$
\ind_F (c) := \frac{\ind_L(c)}{m(c)}.
$$

The closed orbits of $\Phi$ project down to closed geodesics on $X$.
To a geodesic $c$ we attach its free homotopy class of closed paths $[c]$ in
$X$.
The set of free homotopy classes of closed paths stands in bijection to the set
of conjugacy classes of the fundamental group $\Ga$.
We get a map $\{ {\rm closed\ orbits\ of\ } \Phi\}\ra\Ga / {\rm conjugation}$.

Fix a finite dimensional unitary representation $(\ph ,V_\ph)$ of $\Ga$.
Via the above map it makes sense to speak of the number $\tr \ph (c)$ for a
closed orbit $c$.
There is also a geometric interpretation of $\tr\ph(c)$: the representation
$\ph$ defines a flat vector bundle $E_\ph$ over $X$ which lifts to a flat
vector bundle $F_\ph$ over $Y$.
Each orbit $c$ lifts via parallel transport to a linear map $\ph(c)_y$ on the
fibre $F_{\ph,y}$ for any point $y$ on the orbit.
Then the trace $\tr\ph(c)_y$ does not depend on the point $y$ and equals the
above $\tr\ph(c)$.

In the following we will sometimes assume the representation $\ph$ to be {\bf
acyclic}, i.e. $H^q(\Ga ,\ph)=0$ for all $q$.
Under these circumstances we will call the bundle $E_\ph$ acyclic as well.

For $a,b$ nonnegative integers consider the space of smooth sections:
$$
\CC_\ph^{a,b} := \Ga^\infty (\wedge^a T_s^* \otimes \wedge^b T_u^*\otimes
F_\ph).
$$
Let the differential $d_1 : \CC_\ph^{a,b} \ra \CC_\ph^{a+1,b}$ be the
projection of the exterior differential given by the flat connection.
Let $d_2 : \CC_\ph^{a,b} \ra \CC_\ph^{a,b+1}$ be the zero differential.
Form the corresponding bicomplex $(\CC_\ph^{*,*} ,d_1+d_2)$ and its total
complex
$$
\CC_\ph^p := \bigoplus_{a+b=p} \CC_\ph^{a,b}.
$$
Let $H^p(\CC_\ph)$ denote the $p$-th cohomology space.
Note that all the spaces $\CC_\ph^p$ carry natural structures of Fr\'echet
spaces.

For any complex $E=E^0\ra\dots\ra E^n$ over some additive category let $\chi(E)
:=\sum_{p=0}^n(-1)^n H^p(E)$ be the corresponding {\bf Euler object}.
Most generally you might interprete it as a virtual object of the category,
i.e. an element of the free abelian group generated by isomorphism classes of
objects.
Since everything we are going to do with the Euler object factors over exact
sequences we can also consider it as an element of the Grothendieck group.

For operators on infinite dimensional spaces we have the notion of a {\bf
regularized determinant} (see \cite{l2det}).
On a Hilbert space it is defined as follows:
Suppose for an operator $A$ that $A^{-s}$ is defined for $\Re(s)>>0$ and is of
trace class.
Suppose further the zeta function $\zeta_A(s) := \tr (A^{-s})$ extends to a
meromorphic function on $\C$, holomorphic at $s=0$.
Then , extending the finite case, the {\bf determinant} of $A$ is defined as
$$
\det (A) := \exp(-\zeta_A'(0)).
$$
Note that by definition the determinant $\det (A)$ comes with a well defined
logarithm: $\log\det(A) = -\zeta_A'(0)$.

We would like to form the sum over all closed orbits of the Fuller indices but
this sum does not converge.
So we define the {\bf regularized sum} as follows: Consider the sum
$$
\zeta_{\Phi ,\ph}(s) := \sum_c \ind_F(c) \tr \ph(c) e^{-sl(c)},
$$
which converges for $\Re (s)>>0$.
It turns out that $\zeta_{\Phi ,\ph}(s)$ extends to $\C$ with logarithmic
singularities and that in case $\ph$ is acyclic, $\zeta_{\Phi ,\ph}(s)$ is
regular at $s=0$. This means that the number
$$
\hat{\sum_{c}} \ind_F(c) \tr \ph(c) := \zeta_{\Phi ,\ph}(0)
$$
is well defined modulo $2\pi i\Z$.

\begin{theorem} \label{maintheorem}
(Lefschetz formula)
The differentials of the complex $\CC_\ph$ have closed range for any $\ph$, so
the cohomology spaces $H^p(\CC_\ph)$ are Hausdorff.
The unit generator $H$ of the flow $\Phi$ acts on $H^p(\CC_\ph)$ and for
acyclic $\ph$ we have
$$
 \sum_{p=0}^{\dim Y -1} (-1)^{p} \log \det (H|H^P(\CC_\ph)) = -\hat{\sum_{c}}
\ind_F(c)\tr\ph(c) \ \ \ \mod\ 2\pi i\Z.
$$
\end{theorem}

The fact that the differentials have closed image in this context is well known
and follows from the decomposition of the spaces of sections into isotypes
under the action of the isometry group.

Let $\Delta_{\ph ,p}$ be the p-th Laplacian for the flat bundle $E_\ph$ over
$X$.
Ray and Singer \cite{RS-antors} defined the {\bf analytic torsion} of $\ph$ to
be
$$
\tau(\ph) := \prod_{p=0}^{\dim X} \det(\Delta_{\ph ,p})^{p(-1)^p}.
$$

D. Fried \cite{Fr-tors} has proven that:
$$
\log\tau(\ph)= -\hat{\sum_{c}} \ind_F(c)\tr\ph(c) \ \ \ \mod\ 2\pi i\Z.
$$

Putting these things together we conclude
$$
\prod_{p=0}^{\dim X} \det(\Delta_{\ph,p})^{p(-1)^p} = \prod_{q=0}^{2\dim X -2}
\det(H | H^q(\CC_\ph))^{(-1)^{p}}.
$$

\section{A determinant formula}
In this section we will consider the {\bf Ruelle zeta function}
$$
R_\ph (s) := \prod_{c\ prime} \det\left( 1-e^{-sl(c)} \ph(c)\right) ,
$$
where the product is taken over all prime closed geodesics (i.e. those with
$\mu(c)=1$).
Note that in the notation of section \ref{sect_main_theo} we have $\log R_\ph
(s) = \zeta_{\Phi ,\ph}(s)$.

Our Lefschetz formula will follow from

\begin{theorem} \label{detformel}
(determinant formula)
The function $R_\ph(s)$ extends to a meromorphic function on $\C$.
For $\ph$ acyclic it is regular at $s=0$ and it satisfies
$$
R_\ph(s) = e^{P(s)}\prod_{p=0}^{2\dim X -2} \det\left( H+s |
H^p(\CC_\ph)\right)^{(-1)^p},
$$
where $P(s)$ is an odd polynomial.
\end{theorem}

Theorem 4.1 in \cite{onsome} represents $R_\ph$ as an alternating product of
Selberg zeta functions $Z_l$.
So our determinant formula will follow from a corresponding one for the
Selberg function which we will state in greater generality in the next section.

\section{A remark on Selberg zeta functions}
Let $\tilde{X}$ be the universal covering of $X$ and let $G$ denote the
identity component of the isometry group of $\tilde{X}$.
Then $G$ is a semisimple Lie group which acts transitively on $\tilde{X}$.
Let $\tilde{Y}:=S\tilde{X}$ be the sphere bundle of $\tilde{X}$ then via the
differential the action of $G$ lifts to $\tilde{Y}$, since $X$ is a rank one
space this action still is transitive.
So let $M$ be the stabilizer of some point in $\tilde{Y}$ and $(\sigma
,V_\sigma)$ a finite dimensional representation of $M$ then $\sigma$ defines a
$G$-homogeneous vector bundle $\tilde{\CF}_\sigma$ and all such arise this way.
The action of the flow $\Phi$ on $\tilde{Y}$ lifts to a $G$-equivariant action
on the vector bundle and hence pushes down to an action of the induced bundle
$\CF_\sigma$ over $Y$.
For a closed orbit $c$ and $(y,t)\in c$ let $\sigma(c)_y$ be the induced
automorphism of the fibre $\CF_{\sigma ,y}$. We will speak of the trace or the
determinant of $\sigma(c)$ since they do not depend on the point $y$.
Likewise for $(y,t)$ in $c$ let $P_{c,y}^s$ and $P_{c,y}^u$ be the stable and
unstable parts of the
linear Poincar\'e map.

Now we are ready to define the {\bf Selberg zeta function} for $\Re (s)>>0$ as
$$
Z_{\sigma ,\ph}(s) := \prod_{c\ prime} \prod_{N\geq 0} \det\left(
1-e^{-sl(c)}\sigma(c)\otimes\ph(c)
\otimes S^N((P_c^u)^{-1})\right) ,
$$
where $S^N$ means the $N$-th symmetric power.
In \cite{BuOl-buch} it is shown that $Z_{\sigma ,\ph}$ admits meromorphic
continuation, that it satisfies a functional equation and a Riemann hypothesis.

Consider the space of smooth sections:
$$
\CC_{\sigma ,\ph}^a := \Ga^\infty (\wedge^aT_s^* \otimes \CF_\sigma \otimes
F_\ph).
$$
Now $F_\ph$ is equipped with a flat connection, $\tilde{\CF}_\sigma$ has a
unique $G$-invariant connection which pushes down to $\CF_\sigma$.
So we have a canonical connection on the tensor product $\CF_\sigma \otimes
F_\ph$ and we define $d : \CC_{\sigma ,\ph}^a \ra \CC_{\sigma ,\ph}^{a+1}$ as
the projection of the exterior differential given by the connection.
This projection can be expressed in terms of Lie algebra cohomology (see
\cite{onsome}) which shows $d^2=0$.
So the spaces $\CC_{\sigma ,\ph}^.$ form a complex.

\begin{theorem}
Assume that the representation $\sigma$ is self-dual, then for the
Selberg zeta function we have
$$
Z_{\sigma ,\ph}(s) = e^{P_\sigma(s-d+1)}\prod_{p=0}^{\dim X-1} \det(H+s |
H^p(\CC_{\sigma ,\ph}))^{(-1)^{p}},
$$
where $P_\sigma$ is an odd polynomial.
\end{theorem}

\prf
We say that a function $f$ on $\C$ is of {\bf determinant type}
if $f(z)=e^{P(z)} \det(A+z)$ for a polynomial $P(z)$ and an operator
$A$.
The determinant formula in \cite{BuOl-buch} shows that $Z_{\sigma ,\ph}$ is of
determinant type.
By Proposition 3.6 in \cite{onsome} we get
$$
Z_{\sigma ,\ph}(s)=e^{P_\sigma(s-d+1)} \det(H+s | \chi(\CC_{\sigma
,\ph}))^{(-1)^{\dim X-1}}
$$
for some polynomial $P_\sigma$.

The absence of discrete series and the vanishing theorem of $\n$-cohomology
\cite{HeSch} show that the space $H^p(\CC_{\sigma ,\ph})$ are finite
dimensional unless $p=\dim X-1$.
Hence the existence of the determinant on $\chi(\CC_{\sigma ,\ph})$ implies the
existence of the determinant on each single $H^p(\CC_{\sigma ,\ph})$.

It remains to show that the polynomial $P_\sigma$ is odd, i.e. that we have
$P_\sigma(-s)=-P_\sigma(s)$.
To this end we consider the asymptotics of $\log Z_{\sigma ,\ph}(s)$ as $s$
tends to $+\infty$.
By definition it is clear that $\log Z_{\sigma ,\ph}(s)$ tends to zero for
$s\ra +\infty$.
By \cite{BuOl-buch} or \cite{onsome} the divisor of $ Z_{\sigma ,\ph}(s)$ lies
in $(-\infty ,2d-2] \cup (d-1+i\R)$ for some $a,b>0$.
Therefore the function $Z_{\sigma ,\ph}(s)$ which is of determinant type is
actually a product of two determinants:
$$
Z_{\sigma ,\ph}(s+d-1) = e^{P_\sigma(s)} Q(s) \det(A+is)\det(A-is),
$$
where $Q$ is a polynomial and the operator $A$ has spectrum contained in
$[0,\infty )$.

{}From \cite{l2det} we take that there are constants $\alpha_\nu$ and $c_\nu$
for $\nu\geq 0$ with $\alpha_\nu$ real and tending to $+\infty$, such that
$$
-\log \det (A \pm is) = \sum_{\alpha_\nu =0} c_\nu (C+\log s \pm
i\frac{\pi}{2})
$$ $$
+ \sum_{\alpha_\nu = -k \in \{ -1,-2,\dots\} } c_\nu \frac{(\mp i)^k}{k!}
\left(\sum_{j=1}^k\rez{j} -\log(s)\mp i\frac{\pi}{2}\right) s^k
$$ $$
+ \sum_{\alpha_\nu \neq  0,-1,\dots } c_\nu (\pm i)^{-\alpha_\nu} \Ga
(\alpha_\nu)s^{-\alpha_\nu} + o(1),
$$
as $s\ra +\infty$.
Recall that the constants $c_\nu$ are the coefficients of the asymptotic
expansion of the heat kernel of the operator $A$.
Recall further that, by changing the polynomial $Q$ the operator $A$ is only
defined up to finite rank operators.
Since $\sigma$ is self dual it lies in the image of the restriction map
$\Rep K \ra \Rep M$, \cite{Det}.
The argumentations in \cite{Det} give that up to a finite rank operator
$A$ is given as $\sqrt{D_\sigma + c_\sigma}$ for a constant $c_\sigma$ and a
virtual elliptic
differential operator $D_\sigma$ of order two on $X$.
(In \cite{Det} we did assume even dimensions but for odd dimensions
the same arguments apply, they even become easier by the lack of discrete
series.)
Since $\dim X$ is odd it is known that in the asymptotics
of the heat kernel of $D_\sigma$ and hence of $A$ there is no constant term.
So we see that in the above asymptotic expansion the first summand vanishes.
It follows that if
$\log Q(s) +\log(A+is) +\log(A-is)$  is asymptotic to a polynomial it must be
an odd one.
\qed

To conclude the proof of Theorem \ref{detformel} and hence of Theorem
\ref{maintheorem} recall from \cite{onsome} that we have
$$
R_\ph(s) = \prod_{l=0}^{d-1} Z_{\wedge^l\n ,\ph}(s+2l)^{(-1)^l}.
$$
Now we have that as $M$-representations $\wedge^l\n \cong \wedge^{d-1-l}\n$
and hence $Z_{\wedge^l\n ,\ph}=Z_{\wedge^{d-1-l}\n ,\ph}$.
Further the $M$-representation $\wedge^l\n$ is self dual.
Writing $P_l=P_{\wedge^l\n}$ we conclude
$$
R_\ph(s) = e^{P(s)}\prod_{p=0}^{2\dim X -2} \det\left( H+s |
H^p(\CC_\ph)\right)^{(-1)^p}
$$
with
\begin{eqnarray*}
P(-s) &=& \sum_{l=0}^{d-1} P_l(-s+2l-d+1)\\
	&=& \sum_{l=0}^{d-1} P_{d-1-l}(-s+2l-d+1)\\
	&=& \sum_{l=0}^{d-1} P_{l}(-s-2l+d-1)\\
	&=& -\sum_{l=0}^{d-1} P_{l}(s+2l-d+1)\\
	&=& -P(s)
\end{eqnarray*}
\qed

\end{document}